\documentclass[12pt,a4paper]{article}
\usepackage[cp866]{inputenc}
\usepackage[russian]{babel}
\usepackage[psamsfonts]{amssymb}
\newcounter{mchapter}
\newcounter{msection}
\newcounter{lemma}[msection]
\newcounter{theo}[msection]

\newcounter{property}[msection]
\newcounter{corollary}[msection]
\newcounter{definition}[mchapter]

\newcommand{\mchapter}[1]{ \stepcounter{mchapter}.
          \noindent
}
\par
\medskip

\newcommand{\lit}{REFERENCES}

\begin{document}
\centerline{\bf ON LINEAR REPRESENTATIONS OF SOME EXTENSIONS}

\vspace{0.5cm}

\centerline{ V. G. BARDAKOV, O. V. BRYUKHANOV}

\vspace{1cm}

\medskip

\begin{center}
\begin{minipage}[t]{12cm}
\small{{\bf Abstract.} Yu.~I.~Merzljakov developed a method
 of splittable coordinates
which helps  to verify the  linearity of some groups, he established some
fundamental results using
this method.

In this paper we  use the method of splittable coordinates and find some sufficient
condition under which
the semi--direct  product of two linear groups is linear.
As consequence we get  linearity of some HNN-extensions of a free group,
linearity of the holomorph of the braid group $B_n$, $n \geq 2,$ and free group $F_2$ and linearity of
some Artin groups. In all cases we construct faithful linear representations in the explicit
form.



}
\end{minipage}
\end{center}

\vspace{0.5cm}

A.~I.~Malcev [\ref{Mal}] formulated the problem of characterization of abstract groups which
has  faithful linear representation over some field
(we will said that these groups are linear groups).
A.~I.~Malcev found the necessary and sufficient conditions under which
the abelian group has a faithful linear representation.
It is known some results on linear representations of nilpotent, polycyclic,
solvable and some other groups (see the survey of Yu. I. Merzljakov [\ref{Mer}]).
A.~Lubotzky [\ref{Lub}] found a characterization of abstract groups which
has  faithful linear representation over field of complex numbers
 $\mathbb{C}$.
But this result is hard to use for the concrete groups.

Another direction of investigation  is the following problem: if a group
$G$ is construc\-ted from  linear groups with the help of some group construction
(free product, free product with amalgamation, extension, HNN--extension,
wreath product and so on), then for which conditions
group $G$ is linear? In this direction  V.~L.~Nisnevic [\ref{Nis}] proved
that a free product of two linear groups is a linear group.
B.~A.~F.~Wehrfritz [\ref{Weh}] found some sufficient conditions under which
a free product with amalgamation of linear groups is a linear group.
Yu.~E.~Vapne [\ref{Vap1}] studied a wreath product of linear groups. Also he proved
sufficient conditions of linearity of some nilpotent product [\ref{Vap}].
O.~V.~Bryukhanov [\ref{Br5}] found the necessary and sufficient conditions under which
the nilpotent product
 $A(n)B$, $n\ge 2$, of finitely generated linear groups $A$ and $B$ has
a faithful linear representation over fields of zero and positive characteristics.

But the question of  existence of faithful linear representations for group extensions
and HNN--extensions is still open. We have few results in this area.
In particular, it is not known is it true that a semi--direct product of two finitely generated
linear groups is linear. R.~T.~Volvachev [\ref{Vol}] solved the question about
faithful linear rep\-re\-sen\-ta\-tions of HNN--extensions a cyclic group.
By  well--known result of A.~I.~Malcev [\ref{Mal}] every semi--direct product
of residually finite groups is residually finite.
On the other hand it is known that an extension of a linear group by a linear group can be non-linear.
As an example we can consider a finitely generated  free solvable group
of class 3. This group is not linear but is the extension
of free metabelian group of infinite rank by finitely generated abelian group.

Yu.~I.~Merzljakov [\ref{MYuI2}], [\ref{MYuI4}] developed a method
(the method of splittable coordinates)
which helps  to prove linearity of some groups. Yu.~I.~Merzljakov established some fundamental results using
this method.

In this paper, we will use the method of splittable coordinates and find some sufficient
condition under which
the semi--direct  product of two linear groups is linear (see Theorem 1). Also we establish
linearity of ${\rm Hol} (B_n)$, $n \geq 3$ and  ${\rm Hol} (F_2)$ over field of characteristic
zero from these follows linearity over field of characteristic zero any splittable extension
of braid group $B_n$, $n \geq 3$, and free group $F_2$ by linear group.

In Section 2 we consider  HNN--extensions
$$
F_{\varphi}(X) = \langle X, t~||~t^{-1} x t = \varphi(t),~~x \in X \rangle
$$
of a free group $F(X)$ with a virtually inner automorphism $\varphi$ of $F(X)$
(recall that an automorphism $\varphi$ is called virtually inner if $\varphi^m$ is inner
automorphism of $F(X)$ for some
positive integer $m$). Note that these groups are studied in the papers
 [\ref{BBV03}], [\ref{BBV04}].
We  will construct a faithful linear
representation of  $F_{\varphi}(X)$ (see Theorem 2).
Also we will find the center  of $F_{\varphi}(X)$.

In the last Section we study 2-generated Artin groups and construct a faithful linear
representations over ring of rational  fractions $\mathbb{Q}_p$ with denominators equal  to
powers of prime number $p$.
C.~C.~Squier [\ref{Squ}] established the linearity of 2-generated Artin groups
but he constructed representation over field with transcendental (over $\mathbb{Q}$) numbers.

\medskip

{\bf Acknowledgements.}
We would like to thank M.~V.~Neshchadim for his remarks
on the first draft of this paper. Special thanks goes to the participants of the seminar
``Evariste Galois'' at Novosibirsk State University for their
kind attention to our work. Authors were supported by RFBR (grant~02--01--01118).

\vspace{1cm}

\centerline{\bf \S~1. Method of splittable coordinates}

\vspace{0.5cm}

We remind that some homomorphism $\rho : G \longrightarrow {\rm GL}_n(K)$ of group $G$ into
 general linear group
${\rm GL}_n(K)$ over field $K$
is called {\it the linear representation}. This representation $\rho $ is called {\it faithful}
if the kernel ${\rm Ker}(\rho )$ is trivial.
A group is called {\it linear} if it has some faithful linear representation.
Yu.~I.~Merzljakov [\ref{MYuI2}], [\ref{MYuI4}] developed a method
(the method of splittable coordinates)
which helps to prove linearity of some groups.
The heart of this method is to assign a faithful action of group on finite dimension
vector space of algebra of functions from group
$G$ to field $K$. Let $A(G,K)$ be the algebra of all
 $K$--valuable  functions on $G$. If $f\in A (G,K),\;x\in
G$ then we define the action of $x$ by the rule
 \[f^x (y)=f (x y),\;\;y\in G.\]
Since
$f^{gh}=(f^g)^h$ and $f^e=f$, where $e$ is the unit of $G, g,h \in
G$ then we get the action of $G$ on the algebra $A (G,K).$


We will call any subset of functions  $\{ t_i \} \subseteq A(G, K)$, $i \in I$, $K$-{\it valuable
coordinates} on  $G$,
if for any  $x \in G\setminus \{ e \}$ there is $t_j \in \{ t_i \}$ such that
$t_j(x) \neq t_j(e)$.
We will say that the $K$--valuable  coordinates $t_1,\ldots, t_d$ on $G$,
are {\it splittable coordinates} if the following equality are true: \[t_{\alpha}
(xy)=\sum_{\beta=1}^l f_{\alpha\beta}(x)h_{\beta} (y),\;\;x,y\in
G,\;\;1\leq \alpha\leq d,\] where $f_{\alpha\beta},h_{\beta}$ are
some $K$--valuable functions on  $G$.

In the case when
 the $K$--valuable  splittable coordinates $t_1,\ldots, t_d$ are defined on $G$  we can see that
a vector space
$V$ generated by $G$--orbit of the set of coordinate functions
$\{t_{\alpha}\}$ has finite dimension  because it lies in linear
$K$--envelope of the finite set $\{h_{\beta}\}$. Let
$v_1,\ldots,v_n$ be a basis of $V$ over $K$. Then for every  $x\in G$
we have \[v^x_i=\sum_{j=1}^n \rho_{ij} (x) v_j,\;\;\rho_{ij} (x)\in
K.\]

The linear representation $\rho$ which is defined by the rule  $x\mapsto
(\rho_{ij} (x)),\;x\in G$ is faithful. Namely, if the identical linear
transformation of $V$ corresponding to some element $x$ from $G$
then
$t^x_{\alpha}=t_{\alpha}$, where  $1\leq \alpha\leq d$. In this case
$t_{\alpha} (x)=t_{\alpha}^x (e)=t_{\alpha} (e)$ if $1\leq
\alpha\leq d$ and hence $x=e$.
The following lemma describes some properties of this representation
(see [\ref{MYuI2}], [\ref{MYuI4}]).

\medskip

{\bf Lemma 1}~ (on splittable coordinates). {\it Let $K$ be a field and $G$ be a group
with $K$--valuable  coordinates $t_1,\ldots,t_d$ and
\[t_{\alpha} (xy)=\sum^l_{\beta=1} f_{\alpha\beta} (x) h_{\beta} (y),
\;\;x,y \in G,\;\;1\leq\alpha\leq d,\] where $f_{\alpha\beta},
h_{\beta}$ are  functions with the values in  $K$.
There exists the isomorphism
$\rho : G\rightarrow {\rm GL}_n (K)$ such that}: 1) $\rho^{-1}$ {\it is linear
on $G^{\rho}$,} 2) {\it if the function  $f_{\alpha\beta}$ is  polynomial on
$G$ (in coordinates $t_{\alpha}$) and $h_{\beta}$ is  polynomial on
subset $G_1\subset G$ then  $\rho$ is  polynomial on $G_1$.}

\medskip

We will say that the map
 $f :X \times Y\mapsto Mat_n (K)$ {\it is splittable over $K$ of length $l$} if for all
its components  $f_{ij}$ the following equality are true  \[f_{ij}
(x,y)=\sum_{\nu=1}^l \varphi_{\nu i j} (x) \psi_{\nu i j} (y),\;\;x\in X,
y\in Y, \;\; 1 \leq i, j \leq n, \] where $\varphi_{\nu i j},\psi_{\nu i j}$ are the functions
with the values in $K$. This yields that the coordinates $\{ t_{\alpha} \}$
are splittable if the functions
$(x, y) \longmapsto t_{\alpha}(x y)$ are splittable.

\medskip
The following statement was formulated
in [\ref{MYuI4}] as an exercise.
We  give the proof of this statement.

\medskip

{\bf Lemma 2.} {\it Let $G$ be an extension of a group
$A$ by a group $B$ with the factors  $f(x,y) \in
A,\;x,y\in B,$  and automorphisms $ \hat{x}\in Aut (A),\;x\in B$. If
 $A$ and $B$ are linear over the field $K$ and the maps
\[(x,y)\mapsto f (x,y),\;(a,x)\mapsto a^{ \hat{x}}\]
are splittable then the group  $G$ is linear over $K$.}

{\bf Proof.}~ Since  $G$ is an extension of  $A$ by $B$ then for every
$g\in G$ there exists a single pair  $(b,a),\;\;b\in B, a\in A$ and the product of these pairs is defined by
\[ g_1 g_2 = (b_1,a_1)\cdot (b_2,a_2)=(b_1 b_2, f(b_1,b_2) a_1^{ \widehat{b_2}} a_2).\]

Since  $A$ and $B$ are linear then there exist the faithful linear representations:
$$
\rho_A : A \longrightarrow {\rm GL}_m(K),~~~\rho_B : B \longrightarrow {\rm GL}_n(K).
$$
We take the following functions as coordinates on $A$ and $B$:
$$
t_{ij}^{(A)} (a) = \left( \rho_A(a) \right)_{ij},~~a \in A,~~1 \leq i, j \leq m,~~~
t_{ij}^{(B)} (b) = \left( \rho_B(b) \right)_{ij},~~b\in B,~~1 \leq i, j \leq m.
$$
It is evident that the set of functions
$$
\{ T_{ij}^{(A)}~\vert~1 \leq i, j \leq m \} \bigcup \{ T_{ij}^{(B)}~\vert~1 \leq i, j \leq n \},
$$
where
$$
T_{ij}^{(A)} ((b, a)) = t_{ij}^{(A)} (a),\;\;T_{ij}^{(B)} ((b, a)) = t_{ij}^{(B)} (b),
$$
are coordinates on $G$.

Let us show that these coordinates
$T_{ij}^{(A)},\;T_{ij}^{(B)}$ are  splittable. For this we calculate
\[T_{ij}^{(A)} (g_1 g_2)=t_{ij}^{(A)} (f (b_1,b_2) a_1 ^{ \widehat{b_2}} a_2)=\sum_{k,l}
t_{ik}^{(A)}(f (b_1,b_2)) t_{kl}^{(A)}(a_1^{ \widehat{b_2}})
t_{lj}^{(A)}(a_2).\]
Since the maps
\[(x,y)\mapsto f (x,y),\;\;(a,x)\mapsto a^{ \hat{x}}\]
are splittable then
$$
t_{ik}^{(A)}(f (b_1,b_2)) = \sum_{\nu} \varphi_{\nu i k} (b_1) \psi_{\nu i k} (b_2),~~~
t_{kl}^{(A)}(a_1^{ \widehat{b_2}}) = \sum_{\mu}\  \overline{\varphi}_{\mu k l} (a_1)
\overline{\psi}_{\mu k l} (b_2).
$$
Therefore,
$$
\sum_{k,l}
t_{ik}^{(A)}(f (b_1,b_2)) t_{kl}^{(A)}(a_1^{ \widehat{b_2}})
t_{lj}^{(A)}(a_2) =
$$
\[= \sum_{k,l} (\sum_{\nu} \varphi_{\nu i k} (b_1) \psi_{\nu i k} (b_2))
(\sum_{\mu}\  \overline{\varphi}_{\mu k l} (a_1)
\overline{\psi}_{\mu k l} (b_2)) t_{lj}^{(A)}(a_2)=\]
\[=\sum_{k,l,\mu,\nu}(\varphi_{\nu i k} (b_1) \overline{\varphi}_{\mu k l}(b_1))
(\psi_{\nu i k} (b_2) \overline{\psi}_{\mu k l}(b_2) t_{lj}^{(A)}
(a_2)).\]
Let
$$
\Phi_{\alpha} (g_1) = \varphi_{\nu i k} (b_1) \overline{\varphi}_{\mu k l}(b_1),~~~
\Psi_{\alpha j} (g_2) = \psi_{\nu i k} (b_2) \overline{\psi}_{\mu k l}(b_2) t_{lj}^{(A)}
(a_2),~~~\alpha=(i,k,l,\mu,\nu).
$$
then we have
\[T_{ij}^{(A)} (g_1 g_2) = \sum_{\alpha} \Phi_{\alpha} (g_1) \Psi_{\alpha j} (g_2).\]
The set of functions $\{\Psi_{\alpha j}\}$ is finite since the indexes
$i,j,k,l,\mu,\nu$ may take only finite number of values.
This implies that the coordinates $T_{ij}^{(A)}$ are splittable.

For the functions $T_{ij}^{(B)}$ we have
\[T_{ij}^{(B)} (g_1  g_2) = t_{ij}^{(B)} (b_1  b_2)=\sum_{k} t_{ik}^{(B)} (b_1) t_{kj}^{(B)} (b_2) =
\sum_{k} T_{ik}^{(B)} (g_1) T_{kj}^{(B)} (g_2).\] Therefore
the coordinates  $T_{ij}^{(A)}, T_{ij}^{(B)}$ are splittable and linearity
of $G$ follows from Lemma 1.

Now let us prove the following theorem

\bigskip

{\bf Theorem 1.} {\it Let $K$ be a field of arbitrary characteristic,
$\Phi \leq {\rm GL}_m(K)$, $G \leq {\rm GL}_n(K)$.
If for every $\varphi \in \Phi$ there is $g_{\varphi}
\in {\rm GL}_n(K)$ such that $\varphi^{-1} g \varphi = g_{\varphi}^{-1} g g_{\varphi}$ for all}
$g \in G$ {\it then the group  $\Phi \leftthreetimes G$ has faithful linear representation by matrices of size
$ \leq m^2 + n^4$ over the field} $K$.

{\bf Proof.} For every element from
$\Phi \leftthreetimes G$ there is  only one pair  $(\varphi, g)$, $\varphi \in \Phi$, $g \in G$
which represents this element.
Then on the group $\Phi \leftthreetimes G$ we can define the coordinate functions
$T_{ij}^{(1)},$ $1 \leq i, j \leq n^2,$ $T_{pq}^{(2)},$ $1 \leq p, q \leq n,$
by setting
$$
T_{ij}^{(1)} (\varphi, g) = t_{ij}^{(1)}(\varphi),~~~
T_{pq}^{(2)} (\varphi, g) = t_{pq}^{(2)}(g),
$$
where $t_{ij}^{(1)},$ $t_{pq}^{(2)}$  are the components of faithful linear representations
of  $\Phi$ and ${\rm GL}_n(K)$ correspondingly,
i.~e. the coordinate functions on $\Phi$ and on ${\rm GL}_n(K)$.

For every automorphism $\varphi \in \Phi$ we take element $g_{\varphi} \in {\rm GL}_n(K)$
such that $g^{\varphi} = g_{\varphi}^{-1} g g_{\varphi}$ for all $g \in G$.
We  define this function by $\tau : \varphi \longmapsto g_{\varphi}$.
Since $(\varphi_1, g_1) (\varphi_2, g_2) = (\varphi_1 \varphi_2, g_1^{\varphi_2} g_2)$ then
$$
T_{ij}^{(1)} ((\varphi_1, g_1) (\varphi_2, g_2)) =
T_{ij}^{(1)} (\varphi_1 \varphi_2, g_1^{\varphi_2} g_2) = t_{ij}^{(1)}(\varphi_1 \varphi_2) =
\sum_{k=1}^{n^2} t_{ik}^{(1)}(\varphi_1) t_{kj}^{(1)}(\varphi_2) =
$$
$$
= \sum_{k=1}^{n^2} T_{ik}^{(1)}(\varphi_1, g_1) T_{kj}^{(1)}(\varphi_2, g_2),
$$
$$
T_{pq}^{(2)} ((\varphi_1, g_1) (\varphi_2, g_2)) =
T_{pq}^{(2)} (\varphi_1 \varphi_2, g_1^{\varphi_2} g_2) = t_{pq}^{(2)}(g_1^{\varphi_2} g_2) =
t_{pq}^{(2)}(g^{-1}_{\varphi_2} g_1 g_{\varphi_2} g_2) =
$$
$$
= \sum_{k_1,k_2, k_3=1}^{n} t_{pk_1}^{(2)}(g^{-1}_{\varphi_2}) t_{k_1k_2}^{(2)}(g_1)
t_{k_2k_3}^{(2)}(g_{\varphi_2}) t_{k_3q}^{(2)}(g_2) =
$$
$$
= \sum_{k_1,k_2, k_3=1}^{n} t_{k_1k_2}^{(2)}(g_1) t_{pk_1}^{(2)}(g^{-1}_{\varphi_2})
t_{k_2k_3}^{(2)}(g_{\varphi_2}) t_{k_3q}^{(2)}(g_2) =
$$
$$
=\sum_{k_1,k_2=1}^{n} t_{k_1k_2}^{(2)}(g_1)
\left( \sum_{k_3=1}^{n} t_{p k_1}^{(2)}(g^{-1}_{\varphi_2})
t_{k_2 k_3}^{(2)}(g_{\varphi_2}) t_{k_3 q}^{(2)}(g_2) \right)=
$$
$$
= \sum_{k_1,k_2=1}^{n} T_{k_1k_2}^{(2)}(\varphi_1, g_1)
H_{pk_1k_2q}(\varphi_2, g_2),
$$
where
$$
H_{pk_1k_2q}(\varphi_2, g_2) =
\sum_{k_3=1}^{n} t_{p k_1}^{(2)}(g^{-1}_{\varphi_2})
t_{k_2 k_3}^{(2)}(g_{\varphi_2}) t_{k_3 q}^{(2)}(g_2).
$$
Since the set of functions $\{ T_{ij}^{(1)}, H_{pk_1k_2q} \}$ is finite and it has
$m^2 + n^4$ functions then by Lemma 1 the group $\Phi \leftthreetimes G$
has a faithful linear representation
by the matrices of size  $ \leq m^2 + n^4$ over the field $K$.
This completes the proof.

From this theorem easily follows

\bigskip

{\bf Corollary 1.} {\it Let $K$ be a field of arbitrary characteristic and $G \leq {\rm GL}_n(K)$.
Then the group ${\rm Int}(G) G$ embeds in ${\rm GL}_{2n^4}(K)$.}

\medskip

{\bf Corollary 2.} {\it The holomorph ${\rm Hol}(B_n)$ of the braid group $B_n$, $n \geq 2$,
has a faithful linear representation over a field of characteristic zero.
}

\medskip

{\bf Proof.} As was proved in [\ref{DG}] $|{\rm Aut}(B_n) : {\rm Int}(B_n)| = 2$.
Hence,  $|{\rm Hol}(B_n) : {\rm Int}(B_n) B_n| = 2$. Since the braid group $B_n$
has a faithful linear representation over a field of characteristic zero
 (see [\ref{Big}],
[\ref{Kra}]) then by Corollary 1, the group ${\rm Int}(B_n)  B_n$
also has a faithful linear representation over the same field.
Since ${\rm Hol}(B_n)$ is a finite extension of
 ${\rm Int}(B_n) B_n$ then the required statement follows.

\medskip

{\bf Corollary 3.} {\it The holomorph ${\rm Hol}(F_2)$ of free  group $F_2$
has a faithful linear representation over a field of characteristic zero.
}

\medskip

{\bf Proof.} Since the center of $F_2$ is trivial there is a isomorphism of
${\rm Hol}(F_2)$ in the direct production ${\rm Aut}(F_2) \times {\rm Aut}(F_2)$. To see it
we map ${\rm Aut}(F_2)$ to diagonal subgroup of this production and $F_2$ to subgroup
$1 \times {\rm Int}(F_2)$. Since in ${\rm Aut}(F_2)$ the following formulas
$\varphi^{-1} \widehat{g} \varphi = \widehat{g^{\varphi}}$, where $\varphi \in {\rm Aut}(F_2),$
$\widehat{g} \in {\rm Int}(F_2),$ $g \in F_2,$ $g^{\varphi}$ is the image of $g$ under the action
$\varphi$, are true then the images of $F_2$ and ${\rm Aut}(F_2)$ in
${\rm Aut}(F_2) \times {\rm Aut}(F_2)$ generate a subgroup which is isomorphic to ${\rm Hol}(F_2)$.

Since ${\rm Aut}(F_2)$ is linear over a field of characteristic zero (see [\ref{Bar}], [\ref{DFG}])
then the required assertion follows.

\vspace{1cm}

\centerline{\bf \S~2. On linearity of some HNN--extensions of free group}

\vspace{0.5cm}

In this Section we  consider HNN--extensions
$$
F_{\varphi}(X) = \langle X, t~\vert \vert ~ t^{-1} x t = \varphi (x), x \in X \rangle
$$
of a free group $F(X)$ by an  automorphism
$\varphi \in {\rm Aut}(F(X))$.
We will prove that if the automorphism $\varphi$ is virtually inner
 then $F_{\varphi}(X)$ is linear.
Note that the conjugacy problem for such groups $F_{\varphi}(X)$ was solved in
[\ref{BBV04}].

The following property holds.

\bigskip

{\bf Theorem 2.} {\it Let $F(X)$ be a free group with the a of free generators
$X$, let $\varphi \in {\rm Aut}~F(X)$ be such that
$\varphi^n \in {\rm Int}~F(X)$ for some positive integer $n$. If
$\sigma : F(X) \longrightarrow {\rm GL}_m(P)$ is a faithful linear representation of $F(X)$
by matrices over a field  $P$ then the representation  $\tau$ of $F_{\varphi}(X)$ defining on
generators by the following manner:
$$
t^{\tau} =
\left( \begin{array}{cccccc}
\mathbf{0} & E_m & \mathbf{0} & \cdots & \mathbf{0} & \mathbf{0} \\
\mathbf{0} & \mathbf{0} & E_m &  \cdots & \mathbf{0} & \mathbf{0} \\
 & & \ddots & & &  \\
 & &  & \ddots & &   \\
\mathbf{0} & \mathbf{0} &  \mathbf{0} & \cdots & \mathbf{0} &  E_m \\
s f^{\sigma} & \mathbf{0} & \mathbf{0} &  \cdots & \mathbf{0} & \mathbf{0}
\end{array} \right),
$$
$$
x^{\tau} =
\left( \begin{array}{cccccc}
x^{\sigma} & \mathbf{0} &   \cdots & \mathbf{0} & \mathbf{0} \\
\mathbf{0} & (\varphi^{-1}(x))^{\sigma} &   \cdots & \mathbf{0} & \mathbf{0} \\
 & & \ddots & &   \\
\mathbf{0} & \mathbf{0} &   \cdots & (\varphi^{-n+2}(x))^{\sigma} &  \mathbf{0} \\
\mathbf{0} & \mathbf{0} &   \cdots & \mathbf{0} & (\varphi^{-n+1}(x))^{\sigma}
\end{array} \right),
$$
where  $f \in F(X)$ is such that $\varphi^{n}(x) = f^{-1} x f$, $x \in X,$ and $s \in P^*$
is an arbitrary element of infinite oder, $\mathbf{0}$ and $E_m$ are nullity and unity
matrices by
$m \times m$ and  $t^{\tau}$ and $x^{\tau}$ are matrices by  $mn \times mn$ is a faithful
linear representation.}

{\bf Proof.} Let
$$
F_{\varphi^n}(X) = \langle F(X), t^n \rangle \leq F_{\varphi}(X).
$$
It is obvious that $\vert F_{\varphi}(X) : F_{\varphi^n}(X) \vert = n$, and
$F_{\varphi^n}(X)$ is isomorphic to $\mathbb{Z} \times F(X)$.
Let us define the representation
$$
\overline{\sigma} : F_{\varphi^n}(X) \longrightarrow
\langle sE_m, F^{\sigma}(X) \rangle \leq {\rm GL}_m(K),
$$
where $(t^n)^{\overline{\sigma}}
= s f^{\sigma}$, $x^{\overline{\sigma}} = x^{\sigma},$ $x \in X$ and $s \in P^*$
is an element of infinite order.

The representation $\tau$ of $F_{\varphi}(X)$ will construct with the known manner
from the representation
$\overline{\sigma}$ of subgroup of finite index
$F_{\varphi^n}(X)$ (see  [\ref{Mal}], [\ref{KargM}, p. 164]).

If $\{1, t, t^2, \ldots , t^{n-1} \}$ is the set of right coset representatives
of $F_{\varphi}(X)$ by subgroup
 $F_{\varphi^n}(X)$ then for $g \in F_{\varphi}(X)$
the image $g^{\tau}$ is the matrix of degree $mn$ which consists of
blocks of degree $m$.
The block on position
 $(i,j)$ is equal to  $(t^{i-1} g t^{1-j})^{\overline{\sigma}}$ if
$t^{i-1} g t^{1-j} \in F_{\varphi^n}(X)$ and equal to a null matrix in the opposite case.

In the matrix $t^{\tau}$ the non-null blocks are the only blocks on the positions
$(i,i+1),$ $i = 1, 2, \ldots , n-1$.   The identical matrix $E_m$ stands
on these positions.
On the position $(n,1)$ stands the block $(t^n)^{\overline{\sigma}} =
s f^{\sigma}$. For $x \in X$ in the matrix
$x^{\tau}$ the non-null blocks are the only diagonal-blocks and the block on the position
$(i,i)$ is equal to $(t^{i-1} x t^{1-i})^{\overline{\sigma}} = (\varphi^{-i+1} (x))^{\sigma}$.
This completes the proof.

The following proposition gives the center of $F_{\varphi}(X)$.

\bigskip

{\bf Proposition 1.} {\it Let
$$
F_{\varphi}(X) = \langle X, t~\vert \vert ~ t^{-1} x t = \varphi (x),~~ x \in X \rangle
$$
be an extension  of a free group $F(X)$ by a virtually inner automorphism
$\varphi \in {\rm Aut}(F(X))$. Then the center $Z(F_{\varphi}(X))$ of $F_{\varphi}(X)$ is
the infinite cyclic group which generated by
 $t^{n} w_0$, where $n$ is a minimal natural number for which
$\varphi^{n} \in {\rm Int}(F(X))$, and $w_0$ is defined from the equality}
$\varphi^{n} (x) = w_0 x w_0^{-1},$ $x \in X$.

{\bf Proof.} It is evident that $t^{n} w_0 \in Z(F_{\varphi^n}(X))$ and $(t^{n} w_0)^k =
t^{n k} w_0^k$.
Let
$h \in Z(F_{\varphi}(X))$. As element from  $F_{\varphi}(X)$ it can be presented in the
form $t^l f$,
where $l \in \mathbb{Z}$, $f \in F(X)$. Since $h$ commutes with every
$g \in F_{\varphi}(X)$ then  $h^{-1} g h = f^{-1} t^{-l} g t^l f = g$ or
$t^{-l} g t^l = f g f^{-1}$. Dividing $l$ with residue on $n$ we get $l = n k + r,$
$0 \leq r < n$.
Then $t^{-l} g t^l = t^{-n k} t^{-r} g t^{r} t^{n k} = f g f^{-1}$ or
$t^{-r} g t^r = t^{n k} f g f^{-1} t^{-n k} = t^{n k} w_0^k w_0^{-k} f g f^{-1} w_0^k w_0^{-k}
t^{-n k} = w_0^{-k} f g f^{-1} w_0^k$. Therefore
$r = 0$ and $w_0^{-k} f \in Z(F(X))$. So
$t^l f = t^{n k} w_0^{k} = (t^{n} w_0)^k \in \langle t^{n} w_0 \rangle $.
The proposition is proven.

\vspace{1cm}

\centerline{\bf \S~3. 2-generated Artin groups}

\vspace{0.5cm}

In [\ref{HW}] there was formulated the question on linearity of Artin groups.
In this Section we will consider 2-generated Artin groups $A(m)$ which
have the following genetic code
$$
A(m) = \langle x, y  \, || \, w_m (x, y) = w_m (y,x )
\rangle ,
$$
where $m \geq 3$  and
$$
w_m (u,v) \, = \,
\left\{
\begin{array}{ll}
(u v)^n , &  \mbox{if} \quad m = 2n, \\
(u v)^n \, u , &   \mbox{if} \quad m = 2n+1 .
\end{array}
\right.
$$

As it was proved in [\ref{BBV04}] each 2-generated Artin group $A(m)$ is
an extension of a free group of finite rank by some virtually inner automorphism.
We will use this result and the results of previous Section  to construct
some linear representations of $A(m)$ over field of arbitrary characteristic.
In particular, we will construct the representation over
a ring  $\mathbb{Q}_p$ for every prime $p$, where $\mathbb{Q}_p = \{
\frac{x}{p^n}~\vert~x \in \mathbb{Z}, n \in \mathbb{N} \bigcup \{ 0 \} \}$.
Note, that Squier [\ref{Squ}] constructed a faithful linear representation of dimension 2
for the group $A(m)$
over some extension of $\mathbb{Q}$ which has transcendental elements.

Let $\sigma = \sigma(\lambda, \mu) : F_n \longrightarrow {\rm GL}_2 (\mathbb{Z}[\lambda, \mu])$,
$\lambda, \mu \in \mathbb{C}$ be the two-parameter representation of a free group
$F_n$ in group of matrices ${\rm GL}_2 (\mathbb{Z}[\lambda, \mu])$,
which is defined on the generators by the rule
$$
x_0^{\sigma} \longmapsto
\left( \begin{array}{cc}
1 & 0 \\
\lambda & 1
\end{array} \right),
$$
$$
x_i^{\sigma} \longrightarrow
\left( \begin{array}{cc}
1 & \mu \\
0 & 1
\end{array} \right)^{-i}
\left( \begin{array}{cc}
1 & 0 \\
\lambda & 1
\end{array} \right)
\left( \begin{array}{cc}
1 & \mu \\
0 & 1
\end{array} \right)^{i},~~~i = 1, 2, \ldots, n-1.
$$
It is evident, that this representation is faithful if the group generated by the matrices
$$
\left( \begin{array}{cc}
1 & 0 \\
\lambda & 1
\end{array} \right),~~~
\left( \begin{array}{cc}
1 & \mu \\
0 & 1
\end{array} \right)
$$
is isomorphic to the free group $F_2$.

\medskip

{\bf Proposition 2.} {\it Let $A(2n),$ $n \geq 2,$ be the 2-generated Artin group.
Then the representation
$$
\rho : A(2n) \longrightarrow {\rm GL}_{2n}(\mathbb{Z}[\lambda, \mu, s^{\pm 1}]),~~~
\lambda, \mu \in \mathbb{C},~~s \in \mathbb{C}^*,
$$
which is defined on the generators by the rule
$$
x^{\rho} = {\rm diag} \left(
\left( \begin{array}{cc}
1 &  0 \\
\lambda & 1
\end{array} \right), \ldots ,
\left( \begin{array}{cc}
1 & 0 \\
\lambda & 1
\end{array} \right) \right),~~~
y^{\rho} =
\left( \begin{array}{cccccc}
\mathbf{0} & A & \mathbf{0} & \cdots & \mathbf{0} & \mathbf{0} \\
\mathbf{0} & \mathbf{0} & A &  \cdots & \mathbf{0} & \mathbf{0} \\
 & & \ddots & & &  \\
 & &  & \ddots & &   \\
\mathbf{0} & \mathbf{0} &  \mathbf{0} & \cdots & \mathbf{0} &  A \\
B & \mathbf{0} & \mathbf{0} &  \cdots & \mathbf{0} & \mathbf{0}
\end{array} \right),
$$
where
$$
A = \left( \begin{array}{cc}
1 & -\mu \\
0 & 1
\end{array} \right),~~~
B = s \left( \begin{array}{cc}
1 & 0 \\
-\lambda & 1
\end{array} \right)
\left( \begin{array}{cc}
1 - \lambda \mu & \mu \\
- \lambda & 1
\end{array} \right)^{n-1}
$$
is a linear representation of $A(2n)$. In addition, if $|\lambda| \geq 2$,
$|\mu| \geq 2$ and $s$ is not a root of unity then this representation is faithful.}

{\bf Proof.} We have
$$
A(2n) = \langle x_0, x_1, \ldots , x_{n-1}, t~\vert \vert ~ t^{-1} x_i t = \varphi (x_i),
~~ i = 0, 1, \ldots , n-1 \rangle,
$$
where the automorphism $\varphi$ is defined  by the rule
$$
\varphi(x_0)  = x_0 x_1 \ldots x_{n-2} x_{n-1}
x^{-1}_{n-2} \ldots x_1^{-1} x_0^{-1},~~~\varphi(x_i)  = x_{i-1},~~~i = 1, \ldots , n-1.
$$
In this case $\varphi^n (x_i) = \Delta x_i \Delta^{-1}$, $\varphi(\Delta) = \Delta,$ where
$\Delta = x_0 x_1 \ldots x_{n-1}$ and $x_0,$ $t$ are corresponding to canonical generators of
Artin group, i. e.
$$
A(2n) = \langle x_0, t~\vert \vert ~ (x_0 t)^{n} = (t x_0)^n \rangle,
$$
(see [\ref{BBV04}]).

By Theorem 2 the following matrices are corresponding to generators
$x_0,$ $x_1,$ $\ldots , x_{n-1}, t$:
$$
x_0^{\tau} = {\rm diag} (x_0^{\sigma}, (\varphi^{-1}(x_0))^{\sigma}, \ldots ,
(\varphi^{-n+1}(x_0))^{\sigma}),
$$
$$
\vdots
 $$
 $$
x_{n-1}^{\tau} = {\rm diag} (x_{n-1}^{\sigma}, (\varphi^{-1}(x_{n-1}))^{\sigma}, \ldots ,
(\varphi^{-n+1}(x_{n-1}))^{\sigma}),
$$
$$
t^{\tau} =
\left( \begin{array}{cccccc}
\mathbf{0} & E_2 & \mathbf{0} & \cdots & \mathbf{0} & \mathbf{0} \\
\mathbf{0} & \mathbf{0} & E_2 &  \cdots & \mathbf{0} & \mathbf{0} \\
 & & \ddots & & &  \\
 & &  & \ddots & &   \\
\mathbf{0} & \mathbf{0} &  \mathbf{0} & \cdots & \mathbf{0} &  E_2 \\
s (\Delta^{-1})^{\sigma} & \mathbf{0} & \mathbf{0} &  \cdots & \mathbf{0} & \mathbf{0}
\end{array} \right),
$$
where all matrices have degree  $2n$, and $s \in \mathbb{C}^*$ is not a root of unity.

If we conjugate $t^{\tau}$ and  $x_0^{\tau}$ by the matrix
$$
u = {\rm diag} \left( E_2,
\left( \begin{array}{cc}
1 & -\mu \\
0 & 1
\end{array} \right), \ldots ,
\left( \begin{array}{cc}
1 & -\mu \\
0 & 1
\end{array} \right)^{n-1}
\right),
$$
then we get correspondingly
$$
u^{-1} t^{\tau} u =
\left( \begin{array}{cccccc}
\mathbf{0} & A & \mathbf{0} & \cdots & \mathbf{0} & \mathbf{0} \\
\mathbf{0} & \mathbf{0} & A &  \cdots & \mathbf{0} & \mathbf{0} \\
 & & \ddots & & &  \\
 & &  & \ddots & &   \\
\mathbf{0} & \mathbf{0} &  \mathbf{0} & \cdots & \mathbf{0} &  A \\
B & \mathbf{0} & \mathbf{0} &  \cdots & \mathbf{0} & \mathbf{0}
\end{array} \right),~~~
u^{-1} x_0^{\tau} u  = {\rm diag} \left(
\left( \begin{array}{cc}
1 &  0 \\
\lambda & 1
\end{array} \right), \ldots ,
\left( \begin{array}{cc}
1 & 0 \\
\lambda & 1
\end{array} \right) \right),
$$
where
$$
A = \left( \begin{array}{cc}
1 & -\mu \\
0 & 1
\end{array} \right),~~~
B = s \left( \begin{array}{cc}
1 & 0 \\
-\lambda & 1
\end{array} \right)
\left( \begin{array}{cc}
1 - \lambda \mu & \mu \\
- \lambda & 1
\end{array} \right)^{n-1}.
$$
Let us define a representation $\rho$ taking $x^{\rho} = u^{-1} x_0^{\tau} u,$
$y^{\rho} = u^{-1} t^{\tau} u$.
The elements  $x=x_0$ and $y=t$ are canonical generators of Artin group $A(2n)$.
If in the representation
$\sigma = \sigma (\lambda, \mu)$ we take $|\lambda| \geq 2$, $|\mu| \geq 2,$ and $s$
is a non-zero complex number not equal to a root of unity, then
 $\rho$ is faithful representation.
If we take $\lambda = \mu =2,$ $s = p$ is a prime number then we get a faithful
representation of $A(2n)$,
$n \geq 2$ by matrices of degree $2n$ over the ring $\mathbb{Q}_p$. This completes the proof.

\medskip

In the case when $m = 2 n +1$, $n \geq 1$ is odd number, the 2-generated Artin group $A(m)$
is isomorphic to
$$
F_{2n+1}(\psi) = \langle x_0, x_1, \ldots , x_{2n-1}, t~\vert \vert~ t^{-1} x_i t = \psi (x_i),
~~ i = 0, 1, \ldots , 2n-1 \rangle,
$$
where the automorphism  $\psi$ is given by
$$
\psi(x_0)  = x_0 x_2 \ldots x_{2n-4} x_{2n-2}
x^{-1}_{2n-1} \ldots x_3^{-1} x_1^{-1},~~~\psi(x_i)  = x_{i-1},~~~i = 1, \ldots , 2n - 1.
$$
In this case $\psi^{2(2n+1)}(x_i) = \Sigma x_i \Sigma^{-1}$, where
$$
\Sigma = x_0 x_2 \ldots x_{2n-2} (x_0 x_1 \ldots x_{2n-1})^{-1} x_1 x_3 \ldots x_{2n-1}, ~~~
\psi(\Sigma) = \Sigma
$$
(see [\ref{BBV04}]). Note that the inverse automorphism $\psi^{-1}$ is given by
$$
\psi^{-i} (x_0) = x_i,~~~i = 1, 2, \ldots, 2n - 1,
$$
$$
\psi^{-i} (x_0) = \Sigma^{-1} \psi^{4n+2-i} (x_0) \Sigma,~~~i = 2n, 2n + 1, \ldots, 4n + 1.
$$
The canonical generators of $A(2n+1)$ are  $x = t$ and $y = x_0 t$,
i.~e.
$$
A(2n + 1) = \langle t, x_0 t~ \vert \vert~ (t x_0 t)^n t = (x_0 t^2)^n x_0 t \rangle.
$$

{\bf Proposition 3.} {\it Let $A(2n+1),$ $n \geq 1$ be the 2-generated Artin group
with odd number of generators.
Then the map
$$
\rho : A(2n+1) \longrightarrow {\rm GL}_{4(2n+1)}(\mathbb{Z}[\lambda, \mu, s^{\pm 1}]),~~~
\lambda, \mu \in \mathbb{C},~~s \in \mathbb{C}^*,
$$
which is defined on the generators
$$
x^{\rho} =
\left( \begin{array}{cccccc}
\mathbf{0} & E_2 & \mathbf{0} & \cdots & \mathbf{0} & \mathbf{0} \\
\mathbf{0} & \mathbf{0} & E_2 &  \cdots & \mathbf{0} & \mathbf{0} \\
 & & \ddots & & &  \\
 & &  & \ddots & &   \\
\mathbf{0} & \mathbf{0} &  \mathbf{0} & \cdots & \mathbf{0} &  E_2 \\
s(\Sigma^{-1})^{\sigma} & \mathbf{0} & \mathbf{0} &  \cdots & \mathbf{0} & \mathbf{0}
\end{array} \right),
$$
$$
y^{\rho} =
\left( \begin{array}{cccccc}
\mathbf{0} & x_0^{\sigma} & \mathbf{0} & \cdots & \mathbf{0} & \mathbf{0} \\
\mathbf{0} & \mathbf{0} & (\psi^{-1}(x_0))^{\sigma} &  \cdots & \mathbf{0} & \mathbf{0} \\
 & & \ddots & & &  \\
 & &  & \ddots & &   \\
\mathbf{0} & \mathbf{0} &  \mathbf{0} & \cdots & \mathbf{0} &  (\psi^{-4n}(x_0))^{\sigma} \\
s \left( \Sigma^{-1} \psi(x_0)  \right)^{\sigma} & \mathbf{0} & \mathbf{0} &  \cdots & \mathbf{0} & \mathbf{0}
\end{array} \right),
$$
where $\sigma = \sigma(\lambda, \mu)$ is the representation of a free 2-generated group,
is the linear representation of $A(2n+1)$. In this case, if $|\lambda| \geq 2$,
$|\mu| \geq 2$ and $s$ is not a root of unity then this representation is faithful.}

{\bf Proof.}
By Theorem 2 we get for elements $x_0, x_1, \ldots , x_{2n-1}, t$ the following matrices
$$
x_0^{\tau} = {\rm diag} (x_0^{\sigma}, (\psi^{-1}(x_0))^{\sigma}, \ldots , (\psi^{-(4n+1)}(x_0))^{\sigma}),
$$
$$
\vdots
$$
$$
x_{2n-1}^{\tau} = {\rm diag} (x_{2n-1}^{\sigma}, (\psi^{-1}(x_{2n-1}))^{\sigma}, \ldots ,
(\psi^{-(4n+1)}(x_{2n-1}))^{\sigma}),
$$
$$
t^{\tau} = \left( \begin{array}{cccccc}
\mathbf{0} & E_2 & \mathbf{0} & \cdots & \mathbf{0} & \mathbf{0} \\
\mathbf{0} & \mathbf{0} & E_2 &  \cdots & \mathbf{0} & \mathbf{0} \\
 & & \ddots & & &  \\
 & &  & \ddots & &   \\
\mathbf{0} & \mathbf{0} &  \mathbf{0} & \cdots & \mathbf{0} &  E_2 \\
s (\Sigma^{-1})^{\sigma} & \mathbf{0} & \mathbf{0} &  \cdots & \mathbf{0} & \mathbf{0}
\end{array} \right),
$$
where all matrices have degree $4n+2$, and $s \in \mathbb{C}^*$ is not a root of unity.
The map $x \longmapsto t^{\tau},$ $y \longmapsto x_0^{\tau} t^{\tau}$  define of desired
representation.

In view of Theorem 2, if the representation  $\sigma$ is faithful then  $\tau$ is
the faithful representation of $F_{2n+1}(\psi).$
In this case, the matrices $t^{\tau}$, $x_0^{\tau} t^{\tau}$ are canonical generators of
$A(2n+1)$. This completes the proof.

\medskip

Note, that for $\lambda = \mu = 2$ and $s = p$ is a prime number, we get a representation
of $A(2n+1)$ by matrices of degree $4n+2$ over the ring $\mathbb{Q}_p$.

Since 2-generated Artin group $A(m)$ is a split extension of a free group of finite rank
by a virtually inner automorphism, then using of Lemma 4 we  find the center of  $A(m)$.

\bigskip

{\bf Corollary.} {\it The center of 2-generated Artin group $A(m)$ is infinite cyclic
group. In this case}

a) {\it if $m = 2 n$ then the center $Z(A(m))$ is generated by} $(x y)^n$,

b) {\it if $m = 2 n + 1$ then the center $Z(A(m))$ } is generated by $((x y)^n x)^2$.

{\bf Proof.} a) We have the isomorphism
$$
A(2n) \simeq \langle x_0, x_1, \ldots , x_{n-1}, t~ \vert \vert~ t^{-1} x_i t = \varphi (x_i),~~~
i = 0, 1, \ldots , n-1 \rangle
$$
and $\varphi^n(x_i) = \Delta x_i \Delta^{-1},$ $\Delta = x_0 x_1 \ldots  x_{n-1}$,
$\varphi(\Delta) = \Delta$.
From the definition of $\varphi$ it follows that $n$ is the least positive power for which
the automorphism $\varphi^n$ is inner. Since for this isomorphism
to $t \in F_{\varphi}(X)$
corresponds the element $y \in A(2n)$, and to $x_0 \in F_{\varphi}(X)$ corresponds the element
 $x \in A(2n)$ then to $(x_0 t)^n$ corresponds the element $(x y)^n$.
To use Proposition 1 we note that
$$
(x_0 t)^n = x_0 x_1 \ldots  x_{n-1} t^n = \Delta t^n  = t^n \Delta.
$$

b) We have the isomorphism
$$
A(2n + 1) \simeq \langle x_0, x_1, \ldots , x_{2n-1}, t~ \vert \vert~ t^{-1} x_i t = \psi (x_i),~~~
i = 0, 1, \ldots , 2n-1 \rangle
$$
and $\psi^{2(2n+1)}(x_i) = \Sigma x_i \Sigma^{-1},$ $\Sigma = x_0 x_2 \ldots x_{2n-2}
(x_0 x_1 \ldots x_{2n-1})^{-1} x_1 x_3 \ldots x_{2n-1},$ $\psi{\Sigma} = \Sigma$.
Since the elements  $t = x$ and $x_0 t = y$ are canonical generators of
 $A(2n + 1)$ then element $(t x_0 t)^n t$
maps to $(x y)^n x$. Further, since $((t x_0 t)^n t)^2 = t^{4n+2} \Sigma $ and
$4 n + 2$ is a least positive power in which the automorphism
 $\psi$ is inner then from proposition 1 follows desired statement.

Note that the center of Artin groups was found in the paper of  E.~Briskorn and K.~Saito \cite{BS}.

We will construct a series of representations for  $A(3)$. This group is isomorphic
to the braid group on 3 strings
$B_3$.

Using Lemma 3 we take
$$
X_0 = \left( \begin{array}{cc}
1 & 0 \\
\lambda & 1
\end{array} \right),~~~
X_1 = \left( \begin{array}{cc}
1 & \mu \\
0 & 1
\end{array} \right),~~~\lambda, \mu \in \mathbb{C}.
$$
Then the following matrices
$$
T = \left( \begin{array}{cccccc}
\mathbf{0} & E_2 & \mathbf{0} & \mathbf{0} & \mathbf{0} & \mathbf{0} \\
\mathbf{0} & \mathbf{0} & E_2 &  \mathbf{0} & \mathbf{0} & \mathbf{0} \\
\mathbf{0} & \mathbf{0} & \mathbf{0} &  E_2 & \mathbf{0} & \mathbf{0} \\
\mathbf{0} & \mathbf{0} & \mathbf{0} &  \mathbf{0} & E_2 & \mathbf{0} \\
\mathbf{0} & \mathbf{0} &  \mathbf{0} & \mathbf{0} & \mathbf{0} &  E_2 \\
s \Sigma^{-1} & \mathbf{0} & \mathbf{0} &  \mathbf{0} & \mathbf{0} & \mathbf{0}
\end{array} \right),
$$
$DT,$ where $D = {\rm diag} (X_0, X_1, \psi^{-2} (X_0), \psi^{-3} (X_0), \psi^{-4} (X_0),
\psi^{-5} (X_0))$ correspond to generators.
If we conjugate these matrices by the matrix
$$
U = {\rm diag} (E_2, E_2, \Sigma^{-1}, \Sigma^{-1}, \Sigma^{-1}, \Sigma^{-1}),
$$
then we get
$$
U^{-1} T U = \left( \begin{array}{cccccc}
\mathbf{0} & E_2 & \mathbf{0} & \mathbf{0} & \mathbf{0} & \mathbf{0} \\
\mathbf{0} & \mathbf{0} & \Sigma^{-1} &  \mathbf{0} & \mathbf{0} & \mathbf{0} \\
\mathbf{0} & \mathbf{0} & \mathbf{0} &  E_2 & \mathbf{0} & \mathbf{0} \\
\mathbf{0} & \mathbf{0} & \mathbf{0} &  \mathbf{0} & E_2 & \mathbf{0} \\
\mathbf{0} & \mathbf{0} &  \mathbf{0} & \mathbf{0} & \mathbf{0} &  E_2 \\
s E_2 & \mathbf{0} & \mathbf{0} &  \mathbf{0} & \mathbf{0} & \mathbf{0}
\end{array} \right),
$$
$$
U^{-1} D T U = \left( \begin{array}{cccccc}
\mathbf{0} & X_0 & \mathbf{0} & \mathbf{0} & \mathbf{0} & \mathbf{0} \\
\mathbf{0} & \mathbf{0} & X_1 \Sigma^{-1} &  \mathbf{0} & \mathbf{0} & \mathbf{0} \\
\mathbf{0} & \mathbf{0} & \mathbf{0} &  \psi^4(X_0) & \mathbf{0} & \mathbf{0} \\
\mathbf{0} & \mathbf{0} & \mathbf{0} &  \mathbf{0} & \psi^3(X_0) & \mathbf{0} \\
\mathbf{0} & \mathbf{0} &  \mathbf{0} & \mathbf{0} & \mathbf{0} &  \psi^2(X_0) \\
s \psi(X_0) & \mathbf{0} & \mathbf{0} &  \mathbf{0} & \mathbf{0} & \mathbf{0}
\end{array} \right),~~~
$$
We take these matrices  as canonical generators $X$ and $Y$ of $A(3)$.
It is easily show that the following relation is true:  $X Y X = Y X Y$.

This representation will be faithful if the group generated by
$X_0$ and $X_1$ is
isomorphic to the free group $F_2$. Using direct calculations we find the matrices
$$
\Sigma^{-1} = \left( \begin{array}{cc}
1 - \lambda \mu + \lambda^2 \mu^2 & - \lambda \mu^2 \\
-\lambda^2 \mu & 1 + \lambda \mu
\end{array} \right),~~~
\psi (X_0) = \left( \begin{array}{cc}
1  & - \mu \\
\lambda & 1 - \lambda \mu
\end{array} \right),
$$
$$
\psi^2 (X_0) = \left( \begin{array}{cc}
1 + \lambda \mu  & - \mu \\
\lambda^2 \mu & 1 - \lambda \mu
\end{array} \right),~~~
\psi^3 (X_0) = \left( \begin{array}{cc}
1 + \lambda \mu - \lambda^2 \mu^2 &  \lambda \mu^2 \\
-\lambda (1 - \lambda \mu)^2  & 1 - \lambda \mu + \lambda^2 \mu^2
\end{array} \right),
$$
$$
\psi^4 (X_0) = \left( \begin{array}{cc}
1 - 2 \lambda^2 \mu^2 & \mu (1 + 2 \lambda \mu) \\
\lambda ( 1 + \lambda \mu- 2\lambda^2 \mu^2)  & 1 -\lambda \mu + 2\lambda^2 \mu^2
\end{array} \right),
$$
where $s \in \mathbb{C}^*$  and is not a root of unity.

If in Theorem 2 for $\overline{\sigma}$ we take the representation
$$
(t^n)^{\overline{\sigma}} = \left( \begin{array}{cc}
T(s) & 0 \\
0 & f^{\sigma}
\end{array} \right),~~~
x^{\overline{\sigma}} = \left( \begin{array}{cc}
E_2 & 0 \\
0 & x^{\sigma}
\end{array} \right),~~~
T(s) = \left( \begin{array}{cc}
1 & s \\
0 & 1
\end{array} \right),~~~
$$
where $\sigma : F(X) \longrightarrow {\rm SL}_2(\mathbb{Z})$, $s \in \mathbb{Z} \setminus \{ 0\}$
then we get the representation of $F_{\varphi}(X)$ in ${\rm SL}_{4n}(\mathbb{Z})$. In particular,
for the braid group $B_3$ which isomorphic to
$A(3)$ we can take $\lambda = \mu = 2,$ $s \in \mathbb{Z} \setminus \{ 0\}$
and  get the representation $\overline{\sigma}$ (see Proposition 3)
which is faithful representation by integer valued matrices of dimension 24.

\vskip 20pt


\centerline{\bf \lit}
\begin{enumerate}
\item\label{Mal}
A. I. Malcev, On isomorphic matrix representations of infinite groups.
(Russian; with English summary),
Mat. Sb., 8, № 3 (1940), 405--422.

\item\label{Mer}
Yu. I. Merzljakov,
Linear groups. (Russian, English)
J. Sov. Math., 1, (1973) 571-593; translation from Itogi Nauki, Ser. Mat., Algebra,
Topologija, Geometrija 1970, 75--110
(1971).

\item\label{Lub}
A. A. Lubotzky,  Group--theoretic characterization of linear groups, J. Algebra, 113 (1988),
207--214.

\item\label{Nis}
V. L. Nisnevic,  \"{U}ber Gruppen die durch Matrizen \"{u}ber einem kommutativen Feld isomorph darstellbar
sind (Russian, German summary), Mat. Sb., 8, № 3 (1940), 395--403.

\item\label{Weh}
B. A. F. Wehrfritz,  Generalized free products of linear groups, Proc. London  Math. Soc.,
27, (1973), 402-424.

\item\label{Vap1}
 Yu. E. Vapne,
On the representability of a direct interlacing of groups by matrices. (Russian, English)
Math. Notes 7, (1970), 110--114; translation from Mat. Zametki 7, (1970), 181--189.

\item\label{Vap}
 Yu. E. Vapne,
Matrix representability of nilpotent products of groups. (Russian)
Math. Notes, 14, № 3 (1973), 383--394.

\item\label{Br5} O. V. Bryukhanov, Matrix representability and structure of groups,
Algebra and model theory -- 4, Novosibirsk State Technical University, Novosibirsk, 2003, 15--29.

\item\label{Vol}
R. T. Volvachev,  Linear representability of the groups $G_{n,k,l} = \langle a, t~||~t^{-1} a^k t = a^l
\rangle$. (Russian. English
summary),
Fundam. Prikl. Mat., 4, № 4 (1998), 1415-1418.

\item\label{MYuI2}
Yu. I. Merzljakov, Matrix representation of automorphisms, extensions, and solvable groups,
 Algebra and Logic, 7,  (1970), 169--192.

\item\label{MYuI4}
Yu. I. Merzljakov, Rational groups (Russian), M.: Nauka, 1980.

\item\label{BBV03}
V.  Bardakov, L. Bokut and A. Vesnin, Twisted conjugacy in free groups and Makanin's question,
Southeast Asian Bulletin of Math., 29, № 2 (2005), 209--226.
(See also arXiv:math.GR/0401349.)

\item\label{BBV04}
V. Bardakov, L. Bokut and A. Vesnin, On the conjugacy problem for cyclic extensions of free
groups,
Comm. Algebra, 33, № 6 (2005), 1979--1996.
(See also arXiv:math.GR/0402060.)

\item\label{Squ}
C. C. Squier,  Matrix representations of Artin groups, Proc. Amer. Math. Soc., 103, № 1 (1988),
49--53.

\item\label{Wehr}
B. A. F. Wehrfritz,  Infinite linear groups, Berlin, 1973.

\item\label{DG}
    J.~L. Dyer, E. K. Grossman, The automorphism groups of the braid groups,  Amer. J.  Math.,
103, № 6 (1981), 1151--1169.

\item\label{Big}
    S.~Bigelow, Braid groups are linear, J. Amer.  Math. Soc.,
14, № 2 (2001), 471--486.

\item\label{Kra}
 D. Krammer, Braid groups are linear, Annals of. Math.,
155, № 1 (2002), 131--156.

\item\label{KargM}
M. I. Kargapolov, Yu. I. Merzljakov,
Fundamentals of the theory of groups. Transl. from the 2nd Russian ed. by Robert G. Burns.
(Russian, English) Graduate Texts in Mathematics. 62. New York-Heidelberg-Berlin:
Springer-Verlag. XVII, 203 p. (1979).

\item\label{Bar}
Bardakov V. G. Linear representations of the group of conjugating
automorphisms and the braid groups of some manifolds, Sibirsk. Mat. Zh., 46 N 1 (2005),
   17--31.

Линейные представления группы сопрягающих автоморфизмов и групп
кос некоторых многообразий, Сиб. матем. журн., 46, № 1 (2005), 17--31.

\item\label{DFG}
 Dyer J. L., Formanek E. and   Grossman E. K. On the linearity of
automorphism groups of free groups, Arch. Math., 38, N 5 (1982), 404--409.

\item\label{HW}
T. Hsu, D. T. Wise,  On linear and residual properties of graph products, Michigan Math. J., 46,
№ 2 (1999), 251--259.

\bibitem{BS}
E. Brieskorn, K. Saito.
Artin-Gruppen und Coxeter-Gruppen. (Russian)
Matematika, Moskva 18, № 6 (1974), 56--79.

\end{enumerate}

\begin{minipage}[t]{9cm}
\small{Sobolev Institute of Mathematics,
Novosibirsk 630090, Russia,
e-mail: bardakov@math.nsc.ru }
\end{minipage}

\end{document}